\def\vers{Mar.~22, 2011, v.3g}
\magnification=1200
\hsize=6.5truein
\vsize=8.9truein
\input amssym.tex
\font\bfont=cmr10 at 14pt
\font\bifont=cmti10 at 14pt
\font\lfont=cmr10 at 12pt
\font\mfont=cmr9
\font\sfont=cmr8
\font\mbfont=cmbx9

\def\1{\hskip1pt}
\def\bs{\bigskip}
\def\ms{\medskip}
\def\sk{\smallskip}
\def\nin{\noindent}
\def\ssb{\raise.2ex\h{${\scriptscriptstyle\bullet}$}}
\def\msum{\h{$\sum$}}
\def\mcap{\h{$\bigcap$}}
\def\mcup{\h{$\bigcup$}}

\def\mopl{\h{$\bigoplus$}}
\def\A{{\cal A}}
\def\a{\alpha}
\def\b{\beta}
\def\C{{\bf C}}
\def\d{\partial}
\def\D{{\cal D}}
\def\h{\hbox}
\def\la{\lambda}
\def\Lc{{\cal L}}
\def\m{{\frak m}}
\def\N{{\bf N}}
\def\Oc{{\cal O}}
\def\P{{\bf P}}
\def\q{\quad}
\def\Q{{\bf Q}}
\def\X{\widetilde{X}}
\def\Z{{\bf Z}}
\def\Gr{{\rm Gr}}
\def\nnc{{\rm nnc}}
\def\red{{\rm red}}
\def\sing{{\rm sing}}
\def\codim{\h{\rm codim}}
\def\mult{\h{\rm mult}}
\def\lf{\lfloor}
\def\rf{\rfloor}
\def\into{\hookrightarrow}
\def\hlf{\hbox{$1\over 2$}}
\h{}
\vskip 1cm
\centerline{\bfont On the local zeta functions and the
{\bifont b}-functions}

\ms
\centerline{\bfont of certain hyperplane arrangements}

\bs
\centerline{\lfont(with Appendix by Willem Veys)}

\bs
\centerline{Nero Budur, Morihiko Saito, and Sergey Yuzvinsky}
\footnote{}{\sfont 2000 Mathematics Subject Classification: 14J17, 32S40, 32S22.}
\footnote{}{\sfont The first author is supported by the NSF grant
DMS-0700360.}
\footnote{}{\sfont The second author is partially supported by
Kakenhi 21540037.}

\bs\ms
{\narrower\nin
{\mbfont Abstract.} {\mfont
Conjectures of J.~Igusa for p-adic local zeta functions and of
J.~Denef and F.~Loeser for topological local zeta functions assert
that (the real part of) the poles of these local zeta functions
are roots of the Bernstein-Sato polynomials (i.e. the b-functions).
We prove these conjectures for certain hyperplane arrangements,
including the case of reduced hyperplane arrangements in
three-dimensional affine space.}
\par}

\bs\bs
\centerline{\bf Introduction}

\bs\nin
Let $K$ be a $p$-adic field, i.e.\ a finite extension of $\Q_p$,
and $\Oc_K$ be the ring of integers of $K$.
We have the norm defined by $|x|_K=q^{-v(x)}$ for $x\in K^*$
where $v(x)\in\Z$ is the valuation (or the order) of $x\in K$ and
$q$ is the cardinality of the residue field $\Oc_K/{\m_K}$
with $\m_K$ the maximal ideal.
For a nonconstant polynomial $f\in K[x_1,\dots,x_n]$, Igusa's
$p$-adic local zeta function (associated with the characteristic
function of $\Oc_K^n\subset K^n$, see [Ig1], [Ig4]) is defined by
the meromorphic continuation of the integral
$$Z^p_f(s)=\int_{\Oc_K^n}|f(x)|_K^sdx\q\h{for}\,\,\,{\rm Re}\,s>0.$$
Here $dx$ denotes the Haar measure on the compact open subgroup
$\Oc_K^n$ of $K^n$, which is the $p$-adic analogue of
the polydisk $\Delta^n$ in $\C^n$.
Note that $Z^p_f(s)$ is closely related to the Poincar\'e series
associated with the numbers of solutions of $f=0$ in
$(\Oc_K/\m_K^i)^n$ for $i>0$ in the case $f\in\Oc_K[x_1,\dots,x_n]$.

On the other hand, the Bernstein-Sato polynomial (i.e.\ the
$b$-function) of a polynomial $f\in K[x_1,\dots,x_n]$ is the monic
polynomial $b_f(s)$ of the least degree satisfying the relation
$$b_f(s)f^s=Pf^{s+1}\,\,\,\,\h{in}\,\,\,\,R_f[s]f^s\,\,\,\,
\h{for some}\,\,\,P\in\D_n[s],$$
where $R_f$ is the localization of $R:=K[x_1,\dots,x_n]$ by $f$
and $\D_n$ is the Weyl algebra which is generated over $K$ by
$x_1,\dots,x_n$ and $\d/\d x_1,\dots,\d/\d x_n$.
Here $K$ may be any field of characteristic $0$, and
$b_f(s)$ is invariant by extensions of the field $K$,
see (2.1) below.
(There is a shift of the variable $s$ by $1$ if one uses the
definition of the Bernstein polynomial in [Be] since $f^s$ is
replaced by $f^{s-1}$ there).
The local $b$-function $b_{f,x}(s)$ is defined by replacing
the Weyl algebra $\D_n$ with $\D_{X,x}$.
Note that for a homogeneous polynomial $f$, we have
$b_f(s)=b_{f,0}(s)$.

A conjecture of J. Igusa [Ig2] asserts the following.

\ms\nin
{\bf Conjecture~(A)}$^p$. The real part of any pole of the $p$-adic
local zeta function $Z^p_f(s)$ is a root of $b_f(s)$.

\ms
Inspired by this conjecture, J. Denef and F. Loeser [DL] defined
the topological local zeta function $Z^{top}_{f,x}(s)$ (see
(1.1.1) below) for a nonconstant polynomial $f$ and
$x\in f^{-1}(0)$ in the case $K=\C$, and conjectured the following.

\ms\nin
{\bf Conjecture~(A)}$^{top}$. Any pole of the topological local
zeta function $Z^{top}_{f,x}(s)$ is a root of $b_{f,x}(s)$.

\ms
There is a weaker version of the conjectures, due to Igusa, and
Denef and Loeser respectively, and called the
{\it monodromy conjecture}, as follows.

\ms\nin
{\bf Conjecture~(B)}$^p$. For any pole $\a$ of the $p$-adic
local zeta function $Z_f(s)$, $e^{2\pi i{\rm Re}(\a)}$
is an eigenvalue of the Milnor monodromy of $f_{\C}$ at some
$x\in f_{\C}^{-1}(0)\subset\C^n$ choosing an embedding $K\into\C$,
where $f_{\C}$ is the image of $f$ in $\C[x_1,\dots,x_n]$.

\ms\nin
{\bf Conjecture~(B)}$^{top}$. For any pole $\a$ of the topological
local zeta function $Z^{top}_{f,x}(s)$, $e^{2\pi i\a}$ is an
eigenvalue of the Milnor monodromy of $f$ at $y\in f^{-1}(0)$
sufficiently near $x$.

\ms
In Conjecture~(B)$^p$, it is enough to consider an embedding
$K^f\into\C$, where $K^f$ is the subfield of $K$ generated
by the coefficients of the linear factors of $f$ so that $D$ is
defined over $K_f$.
Originally Conjecture~(A)$^p$ and (B)$^p$ are stated for a
polynomial $f\in F[x_1,\dots,x_n]$ with $F$ a number field and $K$
the completion of $F$ at a prime of $F$ (except possibly for a
finite number of primes).
In the hyperplane arrangement case, however, this assumption does
not seem to be essential since Conjecture~(B)$^p$ is already
proved by [BMT] and Conjecture~(A)$^p$ is reduced to
Conjecture~(C) below.

By Conjecture~(A) we will mean Conjecture~(A)$^p$ or
Conjecture~(A)$^{top}$ depending on whether $K$ is the $p$-adic
or complex number field, and similarly for Conjecture~(B).
Note that the eigenvalues of the Milnor monodromies in
Conjecture~(B) can be defined in a purely algebraic way using
the $V$-filtration of Kashiwara [Ka2] and Malgrange [Ma2] on the
$\D_n[s]$-module $R_f[s]f^s$,
and the union of the eigenvalues of the Milnor monodromies for
$x\in f_{\C}^{-1}(0)\subset\C^n$ is independent of the choice of
the embedding $K^f\into\C$, see (2.1) below.
Moreover we have the following.

\ms\nin
{\bf Proposition~1.} {\it Let $K$ be a subfield of $\C$, and
$f\in K[x_1,\dots,x_n]$.

\sk\nin
$(i)$ A complex number $\la\in\C$ is an eigenvalue of the Milnor
monodromy of $f_{\C}$ at some $x\in f_{\C}^{-1}(0) \subset\C^n$
if and only if there is a root $\a$ of $b_f(s)$ such that
$\la=e^{-2\pi i\a}$.

\sk\nin
$(ii)$ If $K=\C$, then for any $x\in f^{-1}(0)$, there is an open
neighborhood $U$ of $x$ in classical topology such that for any
open neighborhood $U'$ of $x$ in $U$, the following two conditions
are equivalent.

\nin
$(a)$ The number $\la$ is an eigenvalue of the Milnor monodromy of
$f$ at some $y\in f^{-1}(0)\cap U'$.

\nin
$(b)$ There is a root $\a$ of $b_{f,x}(s)$ such that
$\la=e^{-2\pi i\a}$.}

\ms
This follows from [Ka1], [Ma2].
By Proposition~1, Conjecture~(B) can be viewed as the modulo
$\Z$ version of Conjecture~(A), and is weaker than the latter.
It is known that Conjectures~(A) and (B) are rather difficult to
prove, see e.g.\ [ACLM1], [ACLM2], [Den], [DL], [Ig3], [Ig4],
[KSZ], [Lo1], [Lo2], [LVa], [LV1], [LV2], [Ro], [VP], [Ve1], [Ve2],
[Ve3], [Ve4].
For a generalization to the ideal case, see [HMY], [VV]
(using [BMS]).

In this paper we prove Conjecture~(A) for certain affine hyperplane
arrangements $D$ in $K^n$.
Let $D_i$ be the irreducible components of $D$, and $m_i$ be the
multiplicity of $D$ along $D_i$.
Let $f$ be a defining equation of $D$.
Set $d:=\deg D=\deg f=\sum_im_i$.
In [BMT], Conjecture~(A) is reduced to the following.

\ms\nin
{\bf Conjecture~(C)}. Let $D$ be an indecomposable essential
central hyperplane arrangement in $\C^n$ with degree $d$.
Then $b_f(-n/d)=0$.

\ms
Here central and essential respectively mean that $0\in D_i$ for
any $i$ and $\dim \bigcap_iD_i=0$.
We say that $D$ is {\it indecomposable} if it is not a union of the
pullbacks of arrangements by the two projections of some
decomposition $\C^n=\C^{n'}\times\C^{n''}$ as a vector space.
Note that the proof of Conjecture~(B) in [BMT] implies that
$-n/d-1$ is a root of $b_f(s)$ in case $-n/d$ is not,
since the roots of $b_f(s)$ are in $(-2,0)$, see [Sa1].

As for the reduction of Conjecture~(A) to Conjecture~(C) we have
more precisely the following.

\ms\nin
{\bf Theorem~1} [BMT].
{\it For an affine hyperplane arrangement $D$ in $K^n$,
Conjecture~$(A)$ holds if Conjecture~$(C)$ for $(D/L)_{\C}$ holds
for every dense edge $L$ of $D$.}

\ms
Here an edge means an intersection of $D_i$, and $D/L$ denotes the
arrangement in $K^n/L$ defined by the $D_i$ containing $L$ and with
the same multiplicity $m_i$, where we may assume $0\in L$ replacing
the origin of $K^n$ if necessary.
We call an edge $L\ne K^n$ {\it dense} if $D/L$ is indecomposable.
If $K$ is a $p$-adic filed, then $(D/L)_{\C}$ denotes the scalar
extension of $D/L$ defined by choosing an embedding $K^f\into\C$
where $K^f\subset K$ is the smallest subfield such that $f$ and all
the $D_i$ are defined over $K^f$.
We have $(D/L)_{\C}=D/L$ in the case $K=\C$.

Theorem~1 is proved by using a resolution of singularities obtained
by blowing up only the proper transforms of the {\it dense} edges
in [STV] (together with Igusa's calculation of candidates for poles
of the $p$-adic zeta functions [Ig1] in the $p$-adic case, see also
(1.1.3) below).
Because of this very special kind of resolution, all the obtained
candidates for poles contribute at least to the monodromy eigenvalues,
and Conjecture~(B) is proved in [BMT] for all the candidates for
poles using the calculation of the Milnor cohomology of hyperplane
arrangements in [CS] (or [Di], Prop.~6.4.6) together with a
result of [STV] on the
relation between indecomposability and nonvanishing of the Euler
characteristic of the projective complement.
This is contrary to the most other cases where lots of cancelations
of apparent poles occur, see [Den], [Lo1], [Ve1], [Ve2], [Ve3]
(and Remark~(1.2) below).
Recently W.~Veys informed us that there are examples of
hyperplane arrangement of degree $d$ in $\C^n$ such that $-n/d$ is
not a pole of $Z_{f,0}^{top}(s)$ in the case $n=3$ with $D$
non-reduced or $n=5$ with $D$ reduced, see Appendix.
These examples imply a negative answer to Question~(Q) in (1.4).
There are no such examples if $n=2$ or $n=3$ with $D$ reduced,
see Propositions~(1.5) and (1.8) below.

In this paper we prove the following.

\ms\nin
{\bf Theorem~2.} {\it Conjecture~$(C)$ holds in the following cases.

\sk\nin
$(i)$ $\{0\}$ is a good dense edge of $D$.

\sk\nin
$(ii)$ $D$ is reduced with $n\le 3$.

\sk\nin
$(iii)$ $D$ is reduced, $(n,d)=1$, and $D_d$ is generic relative to
the other $D_j\,(j\ne d)$.}

\ms
Here $L$ is called a {\it good dense edge} if for any dense
edges $L'\supset L$, we have
$$n(L)/d(L)\le n(L')/d(L'),$$
where $d(L)=\mult_LD=\msum_{D_i\supset L}m_i$ and
$n(L)=\codim\,L$.
We say that $D_d$ is {\it generic relative to the other}
$D_j\,(j\ne d)$ if any nonzero intersection of $D_j\,(j\ne d)$ is
not contained in $D_d$, see [FT], Example~4.5.

In case (i), Theorem~2 follows from Teitler's refinement [Te]
of Musta\c{t}\v{a}'s formula [Mu] for multiplier ideals using only
dense edges, together with a well-known relation between the jumping
coefficients and the roots of $b_f(s)$, see [ELSV].
In case (ii) or (iii), we use a generalization of Malgrange's
formula for the roots of $b_f(s)$ in the isolated singularity
case (see [Sa1], [Sa2]) reducing the assertion to a certain
combinatorial problem which can be solved under condition (ii) or
(iii), where we need a result from [FT] in case (iii).

Combining Theorems~1 and 2, we get

\ms\nin
{\bf Theorem~3.} {\it For an affine hyperplane arrangement $D$ in
$K^n$, conjecture~$(A)$ holds if for every dense edge $L$ of $D$,
one of the three conditions in Theorem~$2$ is satisfied for
$(D/L)_{\C}$.
In particular, Conjecture~$(A)$ holds in the following cases.

\sk\nin
$(i)$ $D$ is of moderate type.

\sk\nin
$(ii)$ $D$ is reduced with $n\le 3$.

\sk\nin
$(iii)$ $D$ is reduced with $n=4$, and for each $0$-dimensional
dense edge $L$ of $D$, either condition~$(ii)$ or $(iii)$ in
Theorem~$2$ is satisfied for $(D/L)_{\C}$.}

\ms
Here $D$ is called {\it of moderate type} if all the dense edges
are good.
Note that in the case~(iii), condition~(ii) in Theorem~2 is
satisfied for $(D/L)_{\C}$ with $L\ne 0$.
It seems quite difficult to generalize the arguments in this paper
to the non-reduced case even for $n=3$, or to the $4$-dimensional
case even for reduced $D$.

We would like to thank W.~Veys for useful comments and especially
for examples in Appendix solving Question~(Q) in (1.4) negatively.

In Section 1 we recall some facts from the theory of local zeta
functions.
In Section 2 we explain how to calculate the $b$-functions of
homogeneous polynomials, and prove Theorem~2 in cases (i) and (iii).
In Section 2 we prove Theorem~2 in case (ii).
In Appendix by W.~Veys, we describe some examples related to
Question~(Q) in (1.4).

\vfill\eject
\centerline{\bf 1. Local zeta functions}

\bs\nin
{\bf 1.1.}
Let $K$ be the complex or $p$-adic number field.
Let $X$ be a complex manifold of dimension $n$ with $f$ a
holomorphic function on $X$ if $K=\C$, and
$X=K^n$ with $f\in K[x_1,\dots,x_n]$ if $K$ is a $p$-adic field.
Set $D=f^{-1}(0)$.
Let $\sigma:(\X,E)\to(X,D)$ be an embedded resolution
with $E_j$ the irreducible components of $E:=\sigma^*D$. Set
$$E_I^{\circ}=\mcap_{i\in I}\,E_j\setminus\mcup_{i\notin I}\,E_j,\q
m_j=\mult_{E_j}\,\sigma^*D,\q
r_j=\mult_{E_j}\det({\rm Jac}(\sigma)).$$
If $K=\C$, the topological local zeta function for $x\in D$
is defined by
$$Z^{top}_{f,x}(s)=\sum_I\chi(E_I^{\circ}\cap
\sigma^{-1}(x))\prod_{j\in I}{1\over m_js+r_j+1},
\leqno(1.1.1)$$
which is independent of the choice of the resolution (see [DL]).
So we get candidates for poles
$$\a_j:=-{r_j+1\over m_j}.
\leqno(1.1.2)$$
Note that each $\a_j$ is not necessarily a pole of
$Z^{top}_{f,x}(s)$ in general.
It is not easy to determine exactly false poles since there are
cancelations of poles in many cases, see [Den], [Lo1], [Ve1], [Ve2],
[Ve3] (and Remark~(1.2) below).
In the hyperplane arrangement case, however, there is a special
kind of resolution by [STV] so that Conjecture~(B) is proved for
the above candidates for poles although it is still unclear whether
they are really poles.

The situation is similar in the $p$-adic case where Igusa's
calculation (see e.g.\ [Ig4], Theorem~8.2.1 or [Den]) implies that
the poles of the local zeta function are among the
complex numbers
$$\a_{j,k}:=-{r_j+1\over m_j}-{2\pi\sqrt{-1}\1 k\over m_j\log q}
\,\,\,(k\in\Z).
\leqno(1.1.3)$$

\ms\nin
{\bf 1.2.~Remark.}
It is known that there are remarkable cancelations of poles
by the summation in the definition (1.1.1). So it is
not easy to eliminate false poles, although the curve case is
rather well understood, see [Den], [Lo1], [Ve1], [Ve2], [Ve3].
(For a relatively simple proof of Conjecture~(B) for $n=2$,
see [Ro].)
It is also known that only a few of the roots of $b_f(s)$ can be
detected by the local zeta function.

\ms\nin
{\bf 1.3.~Proposition.} {\it
Let $D$ be a hyperplane arrangement defined by a polynomial $f$.
Then the topological local zeta function $Z_{f,x}^{top}(s)$ is a
combinatorial invariant.}

\ms\nin
{\it Proof.}
By the definition of $Z_{f,x}^{top}(s)$ in (1.1.1), we may assume
$D$ is central, $x=0$.
We have to calculate the Euler characteristic of each open stratum
of a stratification of $\sigma^{-1}(0)$ which is induced from the
canonical stratification of a divisor with normal crossings.
In this case $\sigma$ is obtained by taking first the blow-up
$X'\to X=\C^n$ along the origin of $\C^n$, and then taking the base
change of an embedded resolution of $(\P^{n-1},Z)$ by the
projection $X'\to\P^{n-1}$ where $Z:=\P(D)$.
The Euler characteristic of an open stratum is calculated from
those of the closed strata contained in the closure of the given
stratum.
So the assertion follows by induction on $n$ using [DP] together
with the embedded resolution of $(\P^{n-1},Z)$ obtained by blowing
up along the proper transforms of {\it all\1} the edges of $Z$.
Indeed, any intersection of the proper transforms of exceptional
divisors can be written as a product of embedded resolutions for
certain induced arrangements, see loc.~cit.\ and [BS], Prop.~2.7
in this case.
(If we blow up along only the proper transforms of dense edges,
we can not apply an inductive arguments since there is a problem
as below:
For two dense edges $L\subset L'$ of $Z\subset\P^{n-1}$, $L$ is
not necessarily a dense edge of the induced arrangement in $L'$.)
This finishes the proof of Proposition~(1.3).

\ms\nin
{\bf 1.4.~Analogue of Conjecture~(C).}
The following question arises naturally:

\ms\nin
{\bf Question~(Q).}
Let $D$ be an indecomposable essential central hyperplane
arrangement in $\C^n$ defined by a polynomial $f$ of degree $d$.
Then, is $-n/d$ a pole of $Z_{f,0}^{top}(s)$?

\ms
We have a positive answer to this question if $n=2$ or $n=3$ and
$D$ is reduced, see Propositions~(1.5) and (1.8) below.
Recently, W.~Veys informed us that the answer is negative in
general, more precisely, if $n=3$ with $D$ non-reduced or
$n=5$ with $D$ reduced, see Appendix.

Assume, for example, $n=2$ and $d=\sum_{i=1}^em_i$ with
$m_i=\mult_{D_i}D$. Then
$$Z_{f,0}^{top}(s)={1\over d\1s+2}\Bigl(2-e+\sum_{i=1}^e
{1\over m_is+1}\Bigr).\leqno(1.4.1)$$
This immediately follows from the definition of the zeta function
since the embedded resolution is obtained by one blow-up and $2-e$
is the Euler characteristic of the open stratum in $\P^1$.
So $-2/d$ is a pole of order 2 if and only if $2m_i=d$ for some $i$.
If $-2/d$ is not a pole of order 2, then the coefficient $C_{-2/d}$
of ${1\over ds+2}$ is given by
$$C_{-2/d}=2-e+\sum_{i=1}^e{d\over d-2m_i}.$$

The next Proposition gives a positive answer to Question~(Q) in
(1.4) for $n=2$ where $D$ may be non-reduced.
This is a special case of [Ve3], Prop.~2.8.

\ms\nin
{\bf 1.5.~Proposition} (W.~Veys [Ve3]).
{\it With the above notation, assume $n=2$.
Then $-2/d$ is a pole of $Z_{f,0}^{top}(s)$.
More precisely, if $-2/d$ is not a pole of order $2$,
then $C_{-2/d}>0$ if $\{0\}$ is a good dense edge of $D$, and
$C_{-2/d}<0$ otherwise.}

\ms\nin
{\it Proof.} See [Ve3], Prop.~2.8.

\ms\nin
{\bf 1.6.~Proposition.} {\it 
Assume $n=3$, and $D$ is reduced.
Let $\nu_m\,(m\ge 2)$ be the number of points of $Z:=\P(D)$ with
multiplicity $m$.
Then
$$Z_{f,0}^{top}(s)={1\over ds+3}\Biggl(\chi(\P^2\setminus Z)+
{\chi(Z\setminus Z^{\sing})\over s+1}+\sum_m
\Bigl(2-m+{m\over s+1}\Bigr){\nu_m\over ms+2}\Biggr).$$
In particular, $-3/d$ is the only candidate for the pole of order
$2$ of $Z_{f,0}^{top}(s)$, and is really a pole of order $2$
if and only if $d/3\in\Z$ and $\nu_{2d/3}\ne 0$.
If $-3/d$ is not a pole of order $2$, the coefficient $C_{-3/d}$ of
${1\over ds+3}$ is given by}
$$C_{-3/d}={9\over d-3}\Bigl(d-1+\sum_{m\ne 2d/3}
{m(m-1)\over 2d-3m}\nu_m\Bigr).$$

\ms\nin
{\it Proof.}
Since the embedded resolution of $(\P^2,Z)$ is obtained by blowing
up along the singular points of $Z$, the first assertion follows
from the definition of $Z_{f,0}^{top}(s)$ using the partition of
the summation over $m\ne 2$ and $m=2$.
This implies the second assertion since the coefficient of the
double pole is given up to a nonzero multiplicative constant by
$$2-m+{md\over d-3}={2a-1\over a-1}\ne 0,\,\,\,\h{where}\,\,\,m=2a
\,\,\,\h{with}\,\,\,a:=d/3\in\Z.$$
For the simple pole case, we have
$$C_{-3/d}=\chi(\P^2\setminus Z)+
{\chi(Z\setminus Z^{\sing})\,d\over d-3}+\sum_{m\ne 2d/3}
\Bigl(2-m+{md\over d-3}\Bigr){\nu_md\over 2d-3m}.$$
Here
$$\eqalign{\chi(\P^2\setminus Z)&=3-2d+\msum_m\,(m-1)\nu_m,\cr
\chi(Z\setminus Z^{\sing})&=2d-\msum_m\,m\1\nu_m.}$$
Indeed, the first equality is reduced to the calculation of
$\chi(Z)$ which is obtained by using the short exact sequence
$0\to\Q_Z\buildrel{\iota}\over\into\mopl_i\Q_{Z_i}\to
{\rm Coker}\,\iota\to 0,$ since the cokernel of $\iota$
is supported on the singular points of $Z$ and its rank at $p$
is $m_p-1$ where $m_p$ is the multiplicity of $Z$ at $p$.

Substituting these, we see that $C_{-3/d}$ is given by
$$3-2d+{2d^2\over d-3}+\sum_{m\ne 2d/3}\nu_m
\Bigl(m-1-{md\over d-3}+{(2-m)d\over 2d-3m}+{md^2\over(d-3)(2d-3m)}
\Bigr).$$
After some calculation this is transformed to
$${9\over d-3}\Bigl(d-1+\sum_{m\ne 2d/3}{m(m-1)\over 2d-3m}\nu_m
\Bigr).$$
(The detail is left to the reader.)
This finishes the proof of Proposition~(1.6).

\ms\nin
{\bf 1.7.~Remark.} A strong form of the conjecture in [DL]
predicts that the multiplicity of each root of the zeta function
is at most that of the $b$-function.
In general, the multiplicity of the root $-1$ of the $b$-function
of a reduced essential central hyperplane arrangement is $n$
(see [Sa2], Th.~1), and this settles the problem for the root $-1$.
However, the problem is rather difficult for the roots with
multiplicity $2$ even in the case $n=3$.
In this case the only such root is $-3/d$ with $d/3\in\N$ and
$\nu_{2d/3}\ne 0$ by Proposition~(1.6), but it is not easy to
calculate the $b$-function. (Indeed, the multiplicity is calculated
only in the case $\nu_m=0$ for $m>3$ in loc.~cit.)

\ms
Using Proposition~(1.6) we get the Proposition below which gives a
positive answer to Question~(Q) in (1.4) if $D$ is reduced and $n=3$.
W.~Veys has informed us that he had verified an analogue of it
for the (finer) motivic or Hodge zeta functions.
(Here `finer' means that the non-vanishing of the pole for these
do not imply that for $Z_{f,0}^{top}(s)$ although the converse is
true.)

\ms\nin
{\bf 1.8.~Proposition.} {\it Let $D$ be an indecomposable essential
central hyperplane arrangement of degree $d$ in $\C^3$.
Assume $D$ is reduced.
Then $-3/d$ is a pole of $Z_{f,0}^{top}(s)$.
More precisely, if $-3/d$ is not a pole of order $2$, then the
coefficient $C_{-3/d}$ of $1\over ds+3$ is strictly positive if
$\{0\}$ is a good dense edge of $D$, i.e.\ if $m<2d/3$ for any $m$
with $\nu_m\ne 0$, and $C_{-3/d}$ is strictly negative otherwise.}

\ms\nin
{\it Proof.}
We may assume $n=3$ since the case $n=2$ is trivial.
We may further assume $\{0\}$ is not a good dense edge of $D$,
since the assertion in the good dense edge case easily follows from
Proposition~(1.6).
We may thus assume $\nu_{m_0}\ne 0$ for some $m_0:=2a+e$ with
$0<e<a:=d/3$ where we do not assume $a\in\Z$.
Since the sum of the multiplicities of any two singular points of
$Z$ is at most $d+1$, we have $\nu_{m_0}=1$ and
$$m\le 3a-m_0+1=a-e+1\q\h{for any}\,\,\,m\ne m_0\,\,\,\h{with}
\,\,\,\nu_m\ne 0.$$
By Proposition~(1.6) the assertion $C_{-3/d}<0$ is equivalent to
$${(2a+e)(2a+e-1)\over e}>3(3a-1)+
\sum_{m\le a-e+1}{m(m-1)\over 2a-m}\nu_m.$$
To show the last inequality, we may replace $m(m-1)\over 2a-m$
with $m(m-1)\over a+e-1$, since ${1\over a+e-1}\ge{1\over 2a-m}$
for $m\le a-e+1$.
Using ${d\choose 2}=\sum_m{m\choose 2}\nu_m$, the assertion
is then reduced to
$${(2a+e)(2a+e-1)\over e}>3(3a-1)+
{3a(3a-1)-(2a+e)(2a+e-1)\over a+e-1},$$
i.e.
$$(2a+e)(2a+e-1)(a+2e-1)-3(3a-1)(2a+e-1)e$$
$$=2(2a+e-1)(a-e)(a-e-1)>0.$$
Here $a>e+1$, i.e.\ $m_0=2a+e<d-1$ since $D$ is indecomposable.
So the assertion is proved.

\vfill\eject
\centerline{\bf 2. Calculation of $b$-functions}

\bs\nin
{\bf 2.1.}
For a nonconstant polynomial $f\in K[x_1,\dots,x_n]$ with
${\rm char}\,K=0$, the $b$-function $b_f(s)$ can be defined to be
the minimal polynomial of the action of $s$ on
$$\D_n[s]f^s/\D_n[s]f^{s+1}.$$
This implies that $b_f(s)$ is invariant by extensions of $K$ and
its roots are rational numbers since the last assertion holds for
$K=\C$ by [Ka1].

Let $i_f:X\into X\times {\bf A}_K^1$ denote the graph embedding of
$f$ where $X={\bf A}_K^n$.
Then via the global section functor, $R_f[s]f^s$ is identified
with the direct image by $i_f$ of the $\D_X$-module
$\Oc_X[{1\over f}]$ in the notation of the introduction.
This is compatible with extensions of $K$.
Moreover, the regular holonomic $\D_X$-module
$\Gr_V^{\a}((i_f)_*\Oc_X[{1\over f}])$ corresponds via the global
section functor to $\Gr_V^{\a}(R_f[s]f^s)$, and via the de Rham
functor to the $\la$-eigenspace of Deligne's nearby cycle sheaf
$\psi_f\C_X$ ([De]) with $\la=e^{-2\pi i\a}$ if $K=\C$,
see [Ka2], [Ma2].

This implies that the union of the eigenvalues of the Milnor
monodromies for $x\in f_{\C}^{-1}(0)\subset\C^n$ is independent
of the choice of an embedding $K\into\C$ since the $\a$ are
rational numbers.

\ms\nin
{\bf 2.2.~$b$-functions of homogeneous polynomials.}
Assume that $X=\C^n$ and $f$ is a homogeneous polynomial.
Let $F_f$ denote the Milnor fiber of $f$,
and $H^{n-1}(F_f,\C)_{\la}$ be the $\la$-eigenspace of the Milnor
cohomology by the action of the monodromy $T$, where $n=\dim X$.
Set
$$\lf\a\rf=\max\{k\in\Z\mid k\le\a\},\q
{\bf e}(\a)=\exp(2\pi i\a)\q\h{for}\,\,\,\a\in\Q.$$
By [Sa1], Th.~2, there is a decreasing filtration $P$ on
$H^{n-1}(F_f,\C)_{\la}$ such that
$$b_f(-\a)=0\q\h{if}\q
\Gr_P^{\lf n-\a\rf}H^{n-1}(F_f,\C)_{{\bf e}(-\a)}\ne 0,
\leqno(2.2.1)$$
where $P$ coincides with $\widetilde{P}$ in loc.~cit.\ since $f$ is
homogeneous.

Set $U:=\P^{n-1}\setminus Z$ with $Z:=f^{-1}(0)\subset\P^{n-1}$.
By [Sa1], Prop.~4.9, the filtration $P$ on $H^{n-1}(F_f,\C)_{\la}$
is induced by the pole order filtration $P$ on the meromorphic
extension $\Lc^{(k)}$ of a local system $L^{(k)}$ of rank one on
$U$ such that
$$H^j(U,L^{(k)})=H^j(F_f,\C)_{\la},
\leqno(2.2.2)$$
where $\la=\exp(-2\pi ik/d)$ with $d=\deg f$.
Here the local system $L^{(k)}$ is defined by the decomposition
$$\pi_*\C_{F_f}=\msum_{k=0}^{d-1}\,L^{(k)},$$
where $\pi$ is the canonical projection from the affine Milnor
fiber $F_f:=f^{-1}(1)\subset\C^n$ onto $U\subset\P^{n-1}$, and the
action of the monodromy is the multiplication by $\exp(-2\pi ik/d)$
on $L^{(k)}$ so that (2.2.2) holds, see [CS] or [Di], Prop.~6.4.6.
Since $\P^{n-1}$ is simply connected, the local system $L^{(k)}$
is determined by the monodromies around the irreducible components
$Z_j$ of $Z$.
These are given by the multiplication by $\exp(2\pi im_jk/d)$ where
$m_j$ is the multiplicity of the divisor $Z$ along $Z_j$.

We can identify locally $\Lc^{(k)}$ with $\Oc_Y(*Z)h^{-k/d}$ as
a $\D_Y$-module if $h$ defines locally $Z\subset Y:=\P^{n-1}$.
Then the pole order filtration $P$ on $\Lc^{(k)}$ is defined by
$$\h{$P_i\Lc^{(k)}=\Oc_Yh^{-{k\over d}-i}\,\,\,$ if $\,\,i\ge 0,\,\,$
and $\,\,\,0\,\,\,$ otherwise.}
\leqno(2.2.3)$$
Note that the residue of the logarithmic connection on $P_i\Lc^{(k)}$
at a general point of $Z_j$ is the multiplication by
$$\bigl(-\h{$k\over d$}-i\bigr)m_j.
\leqno(2.2.4)$$

The filtration $P^i=P_{-i}$ on
$H^{n-1}(U,L^{(k)})=H^{n-1}(F_f,\C)_{\la}$ is induced by
$P_{n-1-i}$ on $\Lc^{(k)}$ using the de Rham complex
$$\Lc^{(k)}\to\Lc^{(k)}\otimes_{\Oc_Y}\Omega_Y^1\to\cdots\to
\Lc^{(k)}\otimes_{\Oc_Y}\Omega_Y^{n-1},$$
since the latter has the filtration $P^i=P_{-i}$ defined by
$$P_{-i}\Lc^{(k)}\to P_{1-i}\Lc^{(k)}\otimes_{\Oc_Y}\Omega_Y^1
\to\cdots\to P_{n-1-i}\Lc^{(k)}\otimes_{\Oc_Y}\Omega_Y^{n-1}.$$
We have also the Hodge filtration $F$ on $\Lc^{(k)}$ such that
$$F_i\Lc^{(k)}\subset P_i\Lc^{(k)},$$
and the Hodge filtration $F$ on
$H^{n-1}(U,L^{(k)})=H^{n-1}(F_f,\C)_{\la}$ is induced by the above
formula with $P$ replaced by $F$.

\ms\nin
{\bf 2.3.~Calculation of the cohomology of $L^{(k)}$.}
From now on, assume $D=f^{-1}(0)$ is a central hyperplane
arrangement in $\C^n$.
Let $D_i\,\,(i=1,\dots,e)$ be the irreducible components of $D$
with multiplicity $m_i$.
Then $Z=\P(D)\subset\P^{n-1}$ and $Z_i=\P(D_i)$.
Let $D^{\nnc}$ denote the smallest subset of $D$ such that
$D\setminus D^{\nnc}$ is a divisor with normal crossings.
Set $Z^{\nnc}=\P(D^{\nnc})\subset\P^{n-1}$.
Note that $d=\deg f=\sum_{i=1}^em_i$.

For $k\in\{0,\dots,d-1\}$ and $I\subset\{1,\dots,e-1\}$ with
$|I|=k-1$, define
$$\eqalign{\a^I_i&=\cases{-m_ik/d
&if $i\notin I\cup\{e\}$,\cr
1-m_ik/d&if $i\in I\cup\{e\}$.}\cr
\a^I_L&=\msum_{D_i\supset L}\,\a^I_i.\cr
\Sigma^I&=\{p\in Z^{\nnc}\setminus Z_e\mid\a^I_p=0\},}
\leqno(2.3.1)$$
where $L$ is an edge of $D$, and we set $\a^I_p:=\a^I_L$ if
$\P(L)=\{p\}$.
(See Remark~3.6$\,$(iii) below for another way of the
definition of the $\a^I_i$.)
Here it should be noted that in order to apply the theory in [ESV]
(and also in [STV]), we must have a regular singular connection on
a {\it trivial} line bundle, i.e. the following condition should be
satisfied:
$$\msum_{i=1}^e\,\a^I_i=0.
\leqno(2.3.2)$$
This is satisfied in this case since $d=\msum_{i=1}^e\,m_i$.
Note also that $\a^I_e$ is used in an essential way for (2.3.4)
below (i.e.\ the condition of [STV]) although it does not appear
in the definition of the connection on the affine space $\C^{n-1}$
which is given below.

For $i\in\{1,\dots,e-1\}$, let $e_i=dg_i/g_i$ with $g_i$ a linear
function defining $Z_i\setminus Z_e$ in
$\P^{n-1}\setminus Z_e\cong\C^{n-1}$.
Set
$$\omega_I:=\msum_{i=1}^{e-1}\,\a^I_ie_i.$$
It defines a connection $\nabla^{\omega_I}$ on $\Oc_U$
(where $U=\P^{n-1}\setminus Z$) such that
$$\nabla^{\omega_I}u=du+u\omega_I\q\h{for}\,\,\,u\in\Oc_U.$$
The corresponding local system is isomorphic to $L^{(k)}$
by comparing their local monodromies as remarked in (2.2).
Consider the de Rham cohomology
$H_{\rm DR}^{\ssb}(U,(\Oc_U,\nabla^{\omega_I}))$, which is
calculated by the complex of rational forms $(\Omega_U^{\ssb}(U),
\nabla^{\omega_I})$ since $U$ is affine. Set
$$\A^p=\msum_{i_1<\cdots<i_p}\,\C e_{i_1}\wedge\cdots\wedge e_{i_p}.$$
Then we have a natural inclusion of complexes
$$\iota_I^{\ssb}:(\A^{\ssb},\omega_I\wedge)\hookrightarrow
(\Omega_U^{\ssb}(U),\nabla^{\omega_I}),
\leqno(2.3.3)$$
where the source is called the Aomoto complex.
Note that we have for $a\in\A^0=\C$
$$\iota_I^p(a\1e_{i_1}\wedge\cdots\wedge e_{i_p})=
\iota_I^0(a)\,e_{i_1}\wedge\cdots\wedge e_{i_p},$$
and the image of the injection $\iota_I^0:\A^0\,(=\C)\hookrightarrow
\Gamma(U,\Oc_U)$ depends on the choice of $I$.
Indeed, ${\rm Im}\,\iota_I^0$ depends on the trivialization of the
line bundle $\Lc^{(k)}$ which is determined by $I$, see the proof of
Theorem~(2.5) and Remark~(2.7)(i) below.

By [ESV], [STV], (2.3.3) is a quasi-isomorphism if the following
condition is satisfied:
$$\h{$\a^I_L\notin\Z_{>0}$ for any nonzero
dense edges $\,L\subset D$.}
\leqno(2.3.4)$$

\ms\nin
{\bf 2.4.~Remark.}
Assume $D$ is reduced (i.e.\ $m_i=1$) and $(k,d)=1$.
Then condition (2.3.4) is satisfied for any $I$ with $|I|=k-1$
since $\a^I_L\notin\Z$ for any nonzero edge $L$.
Moreover, this assumption implies that $\psi_{f,\la}\C_X$, the
nearby cycle sheaf with eigenvalue $\la:=\exp(-2\pi ik/d)$,
is supported at the origin.
(Indeed, in case the last assertion is not true, there is
$d'\in(0,d)$ and $k'\in\N$ such that $k/d=k'/d'$.
This follows from the calculation of the Milnor cohomology in (2.2)
to $x\in D\setminus\{0\}$.
Here the degree $d'$ of the defining equation of $D$ at
$x\in D\setminus\{0\}$ becomes strictly smaller.
But this contradicts the assumption $(k,d)=1$.)
The above assertion implies further the vanishing of the lower Milnor
cohomology $H^j(F_f,\C)_{\la}$ for $j<n-1$,
since the nearby cycle sheaf $\psi_{f,\la}\C_X$ is a perverse sheaf
up to the shift of complex by $n-1$.
If moreover $D$ is indecomposable, then we get the nonvanishing of
the highest Milnor cohomology $H^{n-1}(F_f,\C)_{\la}$ by (2.2.2),
since the indecomposability is equivalent to the nonvanishing of
the Euler characteristic $\chi(U)$, see [STV].

\ms
Note that Theorem~4.2(e) in [Sa2] remains valid in the non-reduced
case as follows.

\ms\nin
{\bf 2.5.~Theorem.}
{\it Let $V(I)'$ be the subspace of $\A^{n-1}$ generated by
$e_J:=e_{j_1}\wedge\cdots\wedge e_{j_{n-1}}$ for any
$J=\{j_1,\dots,j_{n-1}\}\subset I$ with $j_1<j_2<\cdots<j_{n-1}$.
Let $V(I)$ be the image of $V(I)'$ in
$H^{n-1}(\A^{\ssb},\omega_I\wedge)$, where $\omega_I$ and
$\a^I=(\a^I_i)$ are as in $(2.3)$.
Assume $V(I)\ne 0$ and $(2.3.4)$ holds.
Then $b_f(-{k\over d})=0$.}

\ms\nin
{\it Proof.}
By (2.2.1) it is enough to show that the image of $e_J$ by the
injection $\iota_I^{n-1}$ in (2.3.3) is contained
$P_0\Lc^{(k)}\otimes_{\Oc_Y}\Omega_Y^{n-1}$ in the notation of (2.2).
Here $P^nH^{n-1}(F_f,\C)_{\la}=0$ since $P_{-1}\Lc^{(k)}=0$.
By definition the image of $a\in\A^0=\C$ by $\iota_I^0$ is a global
section $v_a$ of a free $\Oc_Y$-submodule $\Lc_I$ of $\Lc^{(k)}$ such
that the residue of the connection at the generic point of $Z_i$ is
the multiplication by $\a_i^I$ in (2.3.1).
Set $Z^{I\cup\{e\}}:=\mcup_{k=1}^nZ_{j_k}$ with $j_n:=e$. Then
$$v_a\otimes e_J\in\Lc_I(Z^{I\cup\{e\}})\otimes\Omega_Y^{n-1},$$
since $e_J\in\Omega_Y^{n-1}(\log Z^{I\cup\{e\}})=
\Omega_Y^{n-1}(Z^{I\cup\{e\}})$.
Thus the assertion is reduced to
$$\Lc_I(Z^{I\cup\{e\}})\subset P_0\Lc^{(k)},$$
and this is shown by comparing (2.2.4) and (2.3.1).
Indeed, the eigenvalue of the residue of the connection on
$\Lc_I(Z^{I\cup\{e\}})$ is shifted by $-1$ at the generic point of
$Z_j$ for $j\in I\cup\{e\}$, but it is not smaller than
$-m_jk/d$ even after this shift by (2.3.1).
So Theorem~(2.5) is proved.

\ms\nin
{\bf 2.6.~Proof of Theorem~2 in cases (i) and (iii).}
In case (i), $n/d$ is a jumping coefficient by Teitler's refinement
[Te] of Musta\c{t}\v{a}'s formula [Mu] for multiplier ideals using
only dense edges.
Hence it is a root of $b_f(s)$ up to a sign by [ELSV].

In case (iii), condition (2.3.4) is satisfied for any $I$ with
$|I|=n-1$ since $k=n$ and $(n,d)=1$, see Remark~(2.4) above.
By [FT], Example~4.5, the highest degree cohomology of the Aomoto
complex $H^{n-1}(\A^{\ssb},\omega_I\wedge)$ has a monomial basis
(independently of $I$) under the genericity condition on $D_d$.
Take a subset
$$I=\{i_1,\dots,i_{n-1}\}\subset\{1,\dots,d-1\},$$
such that the corresponding form
$e_I=e_{i_1}\wedge\cdots\wedge e_{i_{n-1}}$ is a member of the
obtained monomial basis.
Since (2.3.4) is satisfied, the image of $e_I$ in the cohomology of
the local system does not vanish.
So the assertion follows from Theorem~(2.5) (i.e.\
[Sa2], Th.~4.2(e)).

\ms\nin
{\bf 2.7.~Remarks.} (i)
In the above argument, the image of $e_I$ by $\iota_I^{n-1}$ is
independent of the choice of $I$ up to a nonzero constant multiple.
Indeed, the injection $\iota_I^0$ in (2.3.3) is defined by using the
trivial line bundle $\Lc_I$ in the proof of Theorem~(2.5) which is
determined by the eigenvalues $\a_i^I$ in (2.3.1).
If we take another $I'\subset\{1,\dots,d-1\}$ with $|I'|=n-1$ and
$e_{I'}\ne 0$, then, using the trivialization given by $\Lc_I$, a
nonzero constant section of $\Lc_{I'}$ is identified with the
rational function $c\1 g_{I'}/g_I$ where
$c\in\C^*$ and $g_I=\prod_{i\in I}g_i$ in the notation of (2.3).
This gives the difference between $\iota_I^j$ and $\iota_{I'}^j$
for any $j$.
So the independence follows since $g_Ie_I=c'g_{I'}e_{I'}$ with
$c'\in\C^*$.

\ms
(ii) We can also identify the image of $e_I$ by $\iota_I^{n-1}$
with an element of the Gauss-Manin system of $f$.
The problem is then closely related to the torsion of the Brieskorn
lattice.

\bs\bs
\centerline{\bf 3. The rank 3 case}

\bs\nin
In this section we assume $n=3$ and give two proofs of the case (ii)
in Theorem~2.
Note that the case $n\le 2$ is well-known.
Indeed, it follows for instance from [Mu], [ELSV].

\ms\nin
{\bf 3.1.~Conditions.} From now on we assume
$$n=k=3.$$
We will write $p\subset i$ if $\{p\}\subset Z_i$, and
set $\a^I_p=\a^I_L$ if $\P(L)=\{p\}$.

In the notation of (2.3.1) we will study the following three
conditions:

\ms\nin
$(a)\q\a^I_p\notin\Z_{>0}\,\,\,\h{for any}\,\,\,p\in Z^{\nnc}=
\P(D^{\nnc})$.

\ms\nin
$(b)\q\exists\,p_0\in\bigl(\mcup_{i\in I}Z_i\bigr)^{\sing}\setminus
Z_e$.

\ms\nin
$(c)\q Z\setminus(Z_e\cup\Sigma^I\cup\{p_0\})$ is connected.

\ms\nin
{\bf 3.2.~Remarks.} (i)
In the case $n=3$, condition $(a)$ coincides with condition (2.3.4)
which implies that the inclusion (2.3.3) is a quasi-isomorphism.
Note that we have always the inequality of the dimensions,
see [LY], Prop.~4.2.

\ms
(ii) For $i,j,k\supset p$, there is a well-known relation
$$e_i\wedge e_j=e_i\wedge e_k-e_j\wedge e_k,
\leqno(3.2.1)$$
which is easily checked by setting $g_i=x$, $g_j=y$ and $g_k=x+y$.
This also follows from the relations of the Orlik-Solomon algebra
which are given by $\partial(e_i\wedge e_j\wedge e_k)$ for
$i,j,k\supset p$, see e.g.\ [OT], p.~60.
As in [BDS], Lemma~1.4, this implies for
$\eta=\sum_{i=1}^{e-1}\b_ie_i$ and
$p\in Z^{\nnc}\setminus Z_e$
$$\h{If $\pi_p(\omega_I\wedge\eta)=0$, then $\a^I_p\b_i=
\b_p\a^I_i$ for any $i\supset p$.}
\leqno(3.2.2)$$
Here $\b_p=\sum_{i\supset p}\b_i$, and $\pi_p(\omega_I\wedge\eta)$
is the $p$-component in the direct sum decomposition in [BDS],
1.3.2
$$H^2(U,\Q)=\mopl_p\,L_p,$$
where $p$ runs over $(Z_{\red})_{\sing}\setminus Z_e$, and
$L_p$ is a vector space of rank $m'_p-1$ with $m'_p$ the
multiplicity of $Z_{\red}$ at $p$.
More precisely $L_p$ has a basis consisting of
$e_i\wedge e_k$ with $i\supset p$ and $i\ne k$ where $k$ is any
fixed member such that $k\supset p$.
This also follows from the definition of the Orlik-Solomon algebra
mentioned after (3.2.1), see e.g.\ [OT], p.~60.

We also get
$$\h{If $p\in(Z_i\cap Z_j)\setminus(Z^{\nnc}\cup Z_e)$, then
$\a^I_i\b_j=\a^I_j\b_i$.}
\leqno(3.2.3)$$
In case $\a^I_i\ne 0$ (i.e.\ $m_i\ne d/3$) for any $i\in I$,
we have by (3.2.2--3)
$$\h{If $\pi_p(\omega_I\wedge\eta)=0$ and $p\notin\Sigma^I$, then
$\b_i/\a^I_i$ is independent of $i\supset p$.}
\leqno(3.2.4)$$

\ms
(iii) Lemma~1.4 in [BDS] or above (3.2.2) is essentially known
to the specialists, see [LY], Lemma~3.1 (and also [Fa], [Li2],
[Yu]).
Here the situation is localized at $p$, i.e.\ the lines not passing
through $p$ are neglected, by using the fact that the
relations of the Orlik-Solomon algebra are of the form
$\partial(e_J)$ for certain $J$ and are compatible with the
decomposition by $p$.

\ms\nin
{\bf 3.3.~Proposition.} {\it With the notation and the assumption
of $(2.3)$, assume $n=k=3$ and there is $I\subset\{1,\dots,e-1\}$
such that $|I|=2$ and conditions~$(a)$, $(b)$ and $(c)$ in $(3.1)$
are satisfied.
Then $b_f(-3/d)=0$ where $f$ is a defining polynomial of $D$.}

\ms\nin
{\it Proof.}
Let $p_0$ be as in condition $(b)$ in (3.1), and assume
the following condition is satisfied:
$$\h{$\pi_p(\omega_I\wedge\eta)=0$ for any $p\ne p_0$.}$$
Then $\eta$ is a multiple of $\omega_I$, i.e. $\b_i/\a^I_i$ is
independent of $i$, see Remark~(3.2)(ii).
So we can apply Theorem~(2.5) (i.e.\ [Sa2], Th.~4.2(e)), and
conclude that $b_f(-3/d)=0$.
This finishes the proof of Proposition~(3.3).

\ms\nin
{\bf 3.4.~One proof of Theorem~2(ii).}
We may assume that $\{0\}$ is not a good dense edge, since we
can apply the case (i) otherwise.
By Proposition~(3.3), it is sufficient to show the following:

\ms\nin
{\bf Assertion.} There is an irreducible component $Z_e$ of $Z$
together with a subset $I\subset \{1,\dots,e-1\}$ such that $|I|=2$
and conditions~$(a)$, $(b)$ and $(c)$ in (3.1) are satisfied
changing the order of $\{1,\dots,e\}$ if necessary.

\ms
Note first that $\a^I_p$ can be an integer only in the case
$d/3\in\Z$.
(Indeed, we have $m_p:=\msum_{i\supset p}\,m_i<d$, and hence
$\a^I_p\equiv 3m_p/d\not\equiv 0$ mod $\Z$ unless $d/3\in\Z$.)
Then the above assertion is shown in the case $d/3\notin\Z$ as
follows.

Since $\a^I_p\notin\Z$ for any $p\in Z^{\nnc}$,
condition~(a) is trivially satisfied and $\Sigma^I=\emptyset$
for any choice of $I$.
Assuming $D$ central and indecomposable, there is $p_0\in Z^{\sing}$
together with $Z_e$ and $I$ satisfying condition~(b).
As for condition~(c), it is not satisfied only in the case there is
$Z_i$ passing through $p_0$ and such that
$Z_i\cap Z_{i'}\subset(Z_i\cap Z_e)\cup\{p_0\}$ for any
$i'\notin\{i,e\}$.
(Otherwise, for any $Z_i$ passing through $p_0$, there is $Z_{i'}$
such that $Z_i\cap Z_{i'}\not\subset Z_e\cup\{p_0\}$.)
In this case every $Z_i$ passes through either $p_0$ or
$Z_i\cap Z_e$.
This implies that $|Z^{\nnc}|=2$ since $D$ is indecomposable.
Then, replacing $Z_e$ with $Z_i$ containing $Z^{\nnc}$, we may take
$p_0$ to be any point of $Z^{\sing}\setminus Z^{\nnc}$ and $I$ is
chosen so that $\{p_0\}=\bigcap_{i\in I}Z_i$.
Thus the assertion is proved in this case.

We may now assume
$$a:=d/3\in\Z.$$
Since $\{0\}$ is not a good dense edge, there is $p_1\in Z^{\nnc}$
with multiplicity $>2a$.
On the other hand, we may assume that there is
$p_2\in Z^{\nnc}$ with $\a^I_{p_2}\in\Z$, i.e.
its multiplicity is divisible by $a$, since otherwise the above
conditions are easily satisfied.
Thus we may assume that there are $p_1,p_2\in Z^{\nnc}$
with multiplicity $2a+1$ and $a$ respectively and hence
$Z^{\nnc}=\{p_1,p_2\}$, since $d=3a$.
So the assertion is proved by the same argument as above.

\ms\nin
{\bf 3.5.~Another proof of Theorem~2(ii).}
It is also possible to prove Theorem~2(ii) by taking
$p_0$ to be the point with multiplicity $m_{p_0}>{2\over 3}d$,
which exists since we may assume that $\{0\}$ is not a good dense
edge as in (3.4).
In this case there is a line $Z_{d-1}$ which is different from the
line at infinity $Z_d$ and does not contain $p_0$ since $D$ is
indecomposable.
Moreover there are at least two lines $Z_1,Z_2$ passing through
$p_0$ such that their intersections with $Z_{d-1}$ are ordinary
double points of $Z$ and furthermore their intersections with
$Z_d$ do not have multiplicity $a$ so that conditions $(a)$ and
$(b)$ in (3.1) are satisfied by setting $I=\{1,2\}$.
Indeed, we have $m_{p_0}>{2\over 3}d$, $d-m_{p_0}\ge 2$, and hence
$d>6$, and moreover the number of lines $Z_i$ such that
$i\supset p_0$ and $Z_i\cap Z_{d-1}$ is an ordinary double point of
$Z$ is at least
$$m_{p_0}-1-(d-2-m_{p_0})>1,$$
since $|\bigcup_{i\supset p_0}Z_i\cap Z_{d-1}|\ge m_{p_0}-1$.
So the condition on the intersection with $Z_{d-1}$ is satisfied.
For the intersection with $Z_d$ we can exclude the case where a
point of $Z$ has multiplicity $a$ since this case has a very
special structure as explained at the end of (3.4) (e.g.\ the
singular points of $Z$ other than this point and $p_0$ are ordinary
double points) so that we can easily choose $Z_1$, $Z_2$ satisfying
the above conditions in this case.

We can then prove Theorem~2(ii) without using Proposition~(3.3)
but using (3.2.1).
Indeed, by Theorem~(2.5) (i.e.\ [Sa2], Th.~4.2(e)),
it is enough to show
$$\h{If $(\msum_i\,\a^I_ie_i)\wedge(\msum_j\,\b_je_j)=c\1e_1
\wedge e_2$ for some $c\in\Q$, then $c=0$.}
\leqno(3.5.1)$$
Under the assumption of (3.5.1) we get by using (3.2.1)
$$\h{$\a^I_{p_0}\b_i=\b_{p_0}\a^I_i\,\,$ if $i\supset p_0$ and
$i>2$.}
\leqno(3.5.2)$$
Here we have $\a^I_{p_0}\ne 0$ since $m_{p_0}>{2\over 3}d$.
So we may assume
$$\h{$\b_i=0\,\,$ if $i\supset p_0$ and $i>2$,}
\leqno(3.5.3)$$
by replacing $\b_i$ with $\b_i-c'\a^I_i$ for any $i$ where
$c':=\b_{p_0}/\a^I_{p_0}$.
(Note that this change of $\b_i$ does not affect the hypothesis of
(3.5.1).)
Since $m_{p_0}>4$, (3.5.2) and (3.5.3) imply
$$\b_1+\b_2=\b_{p_0}=0.$$
On the other hand, by (3.2.3) applied to the intersections of
$Z_1,Z_2$ with $Z_{d-1}$, we get
$$\b_1/\a^I_1=\b_{d-1}/\a^I_{d-1}=\b_2/\a^I_2,$$
where $\a^I_1=\a^I_2\ne 0$ and $\a^I_{d-1}\ne 0$ since $Z$ is
reduced.
So $\b_1=\b_2=0$, and (3.5.1) follows.

\ms\nin
{\bf 3.6.~Remarks.}
(i) It does not seem easy to generalize the above arguments to the
non-reduced case.
If $p_0$ is taken to be the point with the highest multiplicity,
there is an example as follows:
Assume $a>6$, and let
$$f=(xy(x-y))^{a-2}(x+y-z)(x+y-2z)(x+2y-2z)(2x+y-2z)z^2.$$
Here $d=3a$, and there does not exist $I$ such that the argument in
(3.5) can be applied if we set $p_0=(0,0,1)$.
Indeed, let $Z_i\,\,(i=1,\dots,8)$ denote the lines defined by the
linear factors of $f$ respecting the order of the factors, where
$e=8$.
Here $Z_e$ must be the line defined by $z=0$ since conditions $(a)$
and $(b)$ in (3.1) cannot be satisfied otherwise.
Then the singular points of $Z\setminus(\{p_0\}\cup Z_e)$ contained
in $Z_1$ or $Z_2$ have all multiplicity $a$, and moreover
$Z_{\red}$ has multiplicity 3 at these points.
So the argument in (3.5) cannot be applied.

\ms
(ii) For a more complicated example, we might consider the
following:
Let $E$ be an elliptic curve in the dual projective space $\P^2$,
and $G$ be the subgroup of torsion points of order three.
This defines a projective hyperplane arrangement in $\P^2$ with
$e=|G|=9$, see e.g.\ [Li].
Let $G_0$ be a subgroup of $G$ with order 3.
Assume $a>6$.
To the lines corresponding to the elements of $G_0$ we give the
multiplicity $a-2$, while the other lines have multiplicity 1.
Then $d=3a$, and $I\cup\{e\}$ should correspond to $G_0+p\subset G$
for some $p\in G$ in order to satisfy condition $(a)$ in (3.1).
(Indeed, if there are $g_1,g_2\in I\cup\{e\}$ such that their
images in $G/G_0$ are different, then there is $g_3\in G$ such
that the images of $g_1,g_2,g_3$ in $G/G_0$ are all different and
moreover $g_1+g_2+g_3=0$.
The last condition is equivalent to the condition that the three
lines corresponding to $g_1,g_2,g_3$ intersect at one point.
Then condition $(a)$ is not satisfied at this point.)
So $p_0$ is contained in $Z_e$, and hence condition $(b)$ cannot be
satisfied.
Thus we cannot prove a generalization of Theorem~2 in this case
by using Theorem~(2.5) (i.e.\ the generalization of [Sa2],
Th.~4.2(e) to the nonreduced case).

\ms
(iii) In order to apply the theory in [ESV] and [STV], we have
to choose the residues $\a_i$ of the connection satisfying the two
conditions (2.3.2) and (2.3.4).
In our case we have $\a_i=n_i-m_ik/d$ with $n_i\in\Z$ by the
monodromy condition, and $\sum_in_i=k$ since $\sum_im_i=d$.
Then, to satisfy (2.3.2), an easy way is to choose a subset $J$ of
$\{1,\dots,e\}$ with $|J|=k$ and set $n_i=1$ for $i\in J$ and
$n_i=0$ otherwise.
Here there are two possibilities depending on whether $e\in J$
or $e\notin J$.
Since $e$ corresponds to the divisor at infinity, this makes some
difference in the calculation of the Aomoto complex which is defined
on the complement affine space $\C^{n-1}$.
In (2.3.1) we considered the former case where $I=J\setminus\{e\}$.
However, it is also possible to consider the latter case where
$I=J$ so that $|I|=k$ instead of $|I|=k-1$, and
$$\a^I_i=\cases{-m_ik/d&if $i\notin I$,\cr 1-m_ik/d&if $i\in I$.}$$
In the latter case, however, it is usually more difficult to satisfy
the three conditions in (3.1).

\ms
(iv) If $n=3$, $d\le 7$ and ${\rm mult}_pZ=3$ for any
$p\in Z^{\nnc}$ in the notation of (2.3), the $b$-function of a
reduced hyperplane arrangement is calculated in [Sa2].

\vfill\eject
\centerline{{\bf Appendix}}

\ms
\centerline{by Willem Veys}

\ms
\centerline{\mfont University of Leuven, Department of Mathematics}

\centerline{\mfont Celestijnenlaan 200 B, B-3001 Leuven (Heverlee),
Belgium}

\bs\nin
This appendix describes some examples solving Question~(Q)
in (1.4) negatively.
I thank the authors of this paper for writing some details.

\ms\nin
{\bf A.1.~Example}.
We first explain an example of a nonreduced hyperplane arrangement
with $n=3$, $d=9$. Let
$$f=xy(x-y)z^2(x-z)^4.$$
This gives a negative answer to Question~(Q) in (1.4).
Indeed, we have $\chi(U)=1$ by using the affine space defined by
$x\ne 0$.
We have $\chi(Z_i^{\circ})=-1$ except for the line defined by
$x=0$, and the Euler characteristic is 0 for the latter.
So we get
$$\eqalign{&Z_{f,0}^{top}(s)={1\over 9s+3}\Biggl(1-{2\over s+1}
-{1\over 2s+1}-{1\over 4s+1}+\Bigl(-1+{3\over s+1}\Bigr)
{1\over 3s+2}\cr
&\quad+\Bigl(-1+{1\over s+1}+{1\over 2s+1}+{1\over 4s+1}\Bigr)
{1\over 7s+2}+{2\over s+1}\Bigl({1\over 2s+1}+{1\over 4s+1}\Bigr)
\Biggr).}$$
Set $\Phi(s)=(9s+3)Z_{f,0}^{top}(s)$.
Since ${1\over 2s+1}+{1\over 4s+1}$ vanishes by substituting
$s=-1/3$, we get
$$\Phi(-1/3)=1-3+\Bigl(-1+{9\over 2}\Bigr)-3\Bigl(-1+{3\over 2}
\Bigr)=0.$$
So the pole at $-1/3$ vanishes.
(It vanishes also for the motivic or Hodge
zeta function.)

Note that the above example does not give a counterexample to
Conjecture~(C).
This is shown by using Theorem~(2.5) and Remark~(3.2)(ii) below.
Here $p=(0:1:0)\in\P^2$, the line at infinity is $\{y=0\}$, and $I$
corresponds to the two lines with multiplicities 2 and 4.
This assertion is also shown by a calculation using the computer
program Asir.

\ms\nin
{\bf A.2.~Example}.
There is an example of a {\it reduced} hyperplane arrangement with
$n=5$ and $d=10$, giving a negative answer to Question~(Q) in (1.4),
and which is defined by a polynomial $f$ as below:
$$f=(x-y)(x-2y)(x-3y)(x-4y)(x-5y)(x+y+z)zuv(u+v+z).$$
In fact, let $Z_1,Z_2,Z_3$ be closed subvarieties of $Y:={\bf P}^4$
defined by
$$Z_1=\{x=y=z=0\},\quad Z_2=\{u=v=z=0\},\quad Z_3=\{x=y=0\}.$$
Let $\rho:Y'\to Y$ be the composition of the blow-up of $Y$ along
$Z_1,Z_2$ and the blow-up along the proper transform of $Z_3$.
This gives an embedded resolution of $(Y,Z)$ where
$Z:=\{f=0\}\subset Y$.
We have a partition $\{S_i\}_{i=0,\dots,3}$ of $Y=\P^4$ defined by
$$S_0=\{z\ne 0\},\quad S_i=Z_i\,\,(i=1,2),\quad S_3=\{z=0\}
\setminus(Z_1\cup Z_2).$$
Consider the pullback of the partition
$$S'_i:=\rho^{-1}(S_i)\,\,(i=0,\dots,3).$$
Let $x',y',u',v'$ be affine coordinates of $S_0$ defined
respectively by ${x\over z}, {y\over z}, {u\over z}, {v\over z}$.
Then
$$S'_0=\widetilde{\bf C}^2_{x',y'}\times{\bf C}^2_{u',v'},\quad
S'_1=\widetilde{\bf P}^2_{x,y,z}\times{\bf P}^1_{u,v},\quad
S'_2={\bf P}^2_{u,v,z}\times{\bf P}^1_{x,y},
$$
where $\widetilde{\bf C}^2_{x',y'}$ and $\widetilde{\bf P}^2_{x,y,z}$
are respectively the blow-up of ${\bf C}^2_{x',y'}$ and
${\bf P}^2_{x,y,z}$ along $(0,0)$ and $(0:0:1)$.
Here the lower indices $_{x,y}$ etc.\ indicate the coordinates.
Note that each $S_i$ is a union of strata of the stratification
associated to the divisor with normal crossings $\rho^{-1}(Z)$.
So we get
$$Z_{f,0}^{top}(s)=\sum_{i=0}^3{\Psi_i(s)\over 10s+5},$$
where $\Psi_i(s)/(10s+5)$ is the sum of the factors of
$Z_{f,0}^{top}(s)$ associated to the strata contained in $S'_i$.
Since the stratification is compatible with the above product
structure, we get
$$\eqalign{\Psi_0(s)&=\Bigl(4-{9\over s+1}+{5\over(s+1)^2}+
\Bigl(-3+{5\over s+1}\Bigr){1\over 5s+2}\Bigr)
\cdot\Bigl(1-{3\over s+1}+{3\over(s+1)^2}\Bigr),\cr
\Psi_1(s)&={1\over 7s+3}\Bigl(4-{13\over s+1}+{11\over(s+1)^2}+
\Bigl(-3+{5\over s+1}\Bigr){1\over 5s+2}\Bigr)
\cdot\Bigl(-1+{3\over s+1}\Bigr),\cr
\Psi_2(s)&={1\over 4s+3}\Bigl(1-{4\over s+1}+{6\over(s+1)^2}\Bigr)
\cdot\Bigl(-4+{6\over s+1}\Bigr),\cr
\Psi_3(s)&=0.}$$
Indeed, let $Z'_0$ be the divisor on ${\bf P}^2_{x,y,z}$ defined by
the product of linear factors of $f$ which are linear
combinations of $x,y,z$, and similarly for $Z''_0$ with $x,y$
replaced by $u,v$.
Then
$$\chi({\bf P}^2\setminus Z'_0)=4,\,\,\,
\chi({\bf P}^2\setminus Z''_0)=1,\,\,\,
\chi(Z'_0\setminus\hbox{Sing}\,Z'_0)=-13,\,\,\,
\chi(Z''_0\setminus\hbox{Sing}\,Z''_0)=-4,$$
and the number of ordinary double points of $Z'_0$ and $Z''_0$
are respectively $11$ and $6$.
The calculation for $\P^1_{u,v}$ and $\P^1_{x,y}$ is similar,
and we get the formulas for $\Psi_1(s)$ and $\Psi_2(s)$ since
the definition of $\Psi_1(s),\Psi_2(s)$ is compatible with the
above product structure using the formula:
$\chi(X_1\times X_2)=\chi(X_1)\cdot \chi(X_2)$ for topological
spaces $X_1,X_2$.
As for the first terms, note that the codimensions of the centers
$Z_1$, $Z_2$ are $3$, and the multiplicities of $f$ at the generic
points of $Z_1$ and $Z_2$ are respectively $7$ and $4$.
The term $\bigl(-3+{5\over s+1}\bigr){1\over 5s+2}$ comes from
the exceptional divisor of the blow-up along the proper transform
of $Z_3$, where the multiplicity of $f$ at the generic point of
$Z_3$ is $5$ and $Z_3$ has codimension $2$.

The argument is similar for $\Psi_0(s)$.
Here the Euler number of the smooth part and the number of ordinary
double points change since the varieties are restricted to (the
blow-up of) the affine space $\C^2$.
The vanishing of $\Psi_3(s)$ follows from the ${\bf C}^*$-action
on $S'_3=S_3$ compatible with the stratification, which is defined
by $\lambda\1(x:y:u:v)=(\lambda\1x:\lambda\1y:u:v)$ for
$\lambda\in\C^*$.

Substituting $s=-\hlf$ to the above formulas, we get
$$\Psi_0\bigl(-\hlf\bigr)=-8\cdot 7,\quad
\Psi_1\bigl(-\hlf\bigr)=-2\cdot 8\cdot 5,\quad
\Psi_2\bigl(-\hlf\bigr)=17\cdot 8,$$
and hence the pole of $Z_{f,0}^{top}(s)$ at $s=-\hlf$ vanishes.
For the moment it is not clear whether $-\hlf$ is a root of $b_f(s)$.

\bs\bs
\centerline{{\bf References}}

\ms
{\mfont\baselineskip=12pt
\item{[ACLM1]}
Artal Bartolo, E., Cassou-Nogu\`es, P., Luengo, I.\ and
Melle Hern\'andez, A., Monodromy conjecture for some surface
singularities, Ann.\ Sci.\ Ecole Norm.\ Sup.\ (4) 35 (2002),
605--640.

\item{[ACLM2]}
Artal Bartolo, E., Cassou-Nogu\`es, P., Luengo, I.\ and
Melle Hern\'andez, A., Quasi-ordinary power series and their zeta
functions, Mem.\ Amer.\ Math.\ Soc.\ 178 (2005), no. 841.

\item{[Be]}
Bernstein, I.N., Analytic continuation of generalized functions with
respect to a parameter, Funk.\ Anal.\ 6 (1972), 26--40.

\item{[BDS]}
Budur, N., Dimca, A.\ and Saito, M.,
First Milnor cohomology of hyperplane arrangements,
in ``Topology of Algebraic Varieties and Singularities",
Contemporary Mathematics 538 (2011), pp.~279--292.

\item{[BMS]}
Budur, N., Musta\c{t}\u{a}, M.\ and Saito, M., Bernstein-Sato
polynomials of arbitrary varieties, Compos. Math. 142 (2006),
779--797.

\item{[BMT]}
Budur, N., Musta\c{t}\u{a}, M.\ and Teitler, Z.,
The monodromy conjecture for hyperplane arrangements,
(arXiv:0906.1991) to appear in Geom.\ Dedicata.

\item{[BS]}
Budur, N.\ and Saito, M., Jumping coefficients and spectrum of a
hyperplane arrangement, Math.\ Ann.\ 347 (2010), 545--579.

\item{[CS]}
Cohen, D.C.\ and Suciu, A., On Milnor fibrations of arrangements,
J.\ London Math.\ Soc.\ 51 (1995), 105--119.

\item{[DP]}
De Concini, C.\ and Procesi, C., Wonderful models of subspace
arrangements, Selecta Math.\ (N.S.) 1 (1995), 459--494.

\item{[De]}
Deligne, P., Le formalisme des cycles \'evanescents, in SGA7 XIII
and XIV, Lect.\ Notes in Math.\ 340, Springer, Berlin, 1973,
pp.\ 82--115 and 116--164.

\item{[Den]}
Denef, J., Report on Igusa's local zeta function,
S\'eminaire Bourbaki no.~741, Ast\'erisque 201/202/203 (1991),
359-- 386.

\item{[DL]}
Denef, J.\ and Loeser, F., Caract\'eristiques d'Euler-Poincar\'e,
fonctions z\'eta locales et modifications analytiques,
J.\ Amer.\ Math.\ Soc.\ 5 (1992), 705--720.

\item{[Di]}
Dimca, A., Sheaves in Topology,
Universitext, Springer, Berlin, 2004.

\item{[ELSV]}
Ein, L., Lazarsfeld, R., Smith, K.E. and Varolin, D., Jumping
coefficients of multiplier ideals, Duke Math. J. 123 (2004),
469--506.

\item{[ESV]}
Esnault, H., Schechtman V.\ and Viehweg, E., Cohomology of local
systems on the complement of hyperplanes, Inv.\ Math.\ 109 (1992),
557--561.

\item{[Fa]}
Falk, M., Arrangements and cohomology, Ann.\ Combin.\ 1 (1997),
135--157.

\item{[FT]}
Falk, M.\ and Terao, H., $\beta${\bf nbc}-bases for cohomology of
local systems on hyperplane complements,
Trans.\ Amer.\ Math.\ Soc.\ 349 (1997), 189--202.

\item{[HMY]}
Howald, J., Musta\c{t}\v{a}, M.\ and Yuen, C., On Igusa zeta
functions of monomial ideals, Proc.\ Amer.\ Math.\ Soc.\ 135 (2007),
3425--3433.

\item{[Ig1]}
Igusa, J.-i., Complex powers and asymptotic expansions, I and II,
J.\ reine angew.\ Math.\ 268/269 (1974), 110--130 and 278/279 (1975),
307--321.

\item{[Ig2]}
Igusa, J.-i., $b$-functions and $p$-adic integrals, in Algebraic
analysis, Vol.\ I, Academic Press, Boston, MA, 1988, pp.~231--241.

\item{[Ig3]}
Igusa, J.-i., Local zeta functions of certain prehomogeneous vector
spaces, Amer.\ J.\ Math.\ 114 (1992), 251--296.

\item{[Ig4]}
Igusa, J.-i., An introduction to the theory of local zeta functions,
AMS/IP Studies in Advanced Mathematics, 14, American Mathematical
Society, Providence, RI; International Press, Cambridge, MA, 2000.

\item{[Ka1]}
Kashiwara, M., $B$-functions and holonomic systems, Inv.\ Math.\
38 (1976/77), 33--53.

\item{[Ka2]}
Kashiwara, M., Vanishing cycle sheaves and holonomic systems of
differential equations, Lect. Notes in Math.\ 1016, Springer,
Berlin, 1983, pp. 134--142.

\item{[KSZ]}
Kimura, T., Sato, F.\ and Zhu, X.-W., On the poles of $p$-adic
complex powers and the $b$-functions of prehomogeneous vector
spaces, Amer.\ J.\ Math.\ 112 (1990), 423--437.

\item{[LVa]}
Lemahieu, A.\ and Van Proeyen, L., Monodromy conjecture for
nondegenerate surface singularities, to appear in Trans. A.M.S.

\item{[LV1]}
Lemahieu, A.\ and Veys, W., On monodromy for a class of surfaces,
C.\ R.\ Math.\ Acad.\ Sci.\ Paris 345 (2007), 633--638.

\item{[LV2]}
Lemahieu, A.\ and Veys, W., Zeta functions and monodromy for
surfaces that are general for a toric idealistic cluster,
Int.\ Math.\ Res.\ Not.\ 2009, 11--62.

\item{[Li]}
Libgober, A., Characteristic varieties of algebraic curves, in
Applications of algebraic geometry to coding theory, physics and
computation (Eilat, 2001), NATO Sci.\ Ser.\ II Math.\ Phys.\ Chem.,
36, Kluwer Acad.\ Publ., Dordrecht, 2001, pp.~215--254.

\item{[LY]}
Libgober, A.\ and Yuzvinsky, S., Cohomology of the Orlik-Solomon
algebras and local systems, Compos.\ Math.\ 121 (2000), 337--361.

\item{[Lo1]}
Loeser, F., Fonctions d'Igusa $p$-adiques et polyn\^omes de
Bernstein, Amer.\ J.\ Math.\ 110 (1988), 1--21.

\item{[Lo2]}
Loeser, F., Fonctions d'Igusa $p$-adiques, polyn\^omes de Bernstein,
et poly\`edres de Newton, J.\ reine angew.\ Math.\ 412 (1990),
75--96.

\item{[Ma1]}
Malgrange, B., Le polyn\^ome de Bernstein d'une
singularit\'e isol\'ee, in Lect. Notes in Math. 459, Springer,
Berlin, 1975, pp. 98--119.

\item{[Ma2]}
Malgrange, B., Polyn\^ome de Bernstein-Sato et cohomologie
\'evanescente, Analysis and topology on singular spaces, II, III
(Luminy, 1981), Ast\'erisque 101--102 (1983), 243--267.

\item{[Mu]}
Musta\c{t}\v{a}, M., Multiplier ideals of hyperplane arrangements,
Trans. Amer. Math. Soc. 358 (2006), 5015--5023.

\item{[OT]}
Orlik, P.\ and Terao, H., Arrangements of hyperplanes, Springer,
Berlin, 1992.

\item{[Ro]}
Rodrigues, B.\ On the monodromy conjecture for curves on normal
surfaces, Math.\ Proc.\ Cambridge Philos.\ Soc.\ 136 (2004),
313--324.

\item{[Sa1]}
Saito, M., Multiplier ideals, $b$-function, and spectrum of a
hypersurface singularity, Compo.\ Math.\ 143 (2007) 1050--1068.

\item{[Sa2]}
Saito, M., Bernstein-Sato polynomials of hyperplane arrangements
(math.AG/0602527).

\item{[STV]}
Schechtman, V., Terao, H. and Varchenko, A.,
Local systems over complements of hyperplanes and the Kac-Kazhdan
conditions for singular vectors,
J.\ Pure Appl.\ Algebra 100 (1995), 93--102.

\item{[Te]}
Teitler, Z., A note on Musta\c{t}\v{a}'s computation of multiplier
ideals of hyperplane arrangements, Proc.\ Amer.\ Math.\ Soc.\ 136
(2008), 1575--1579.

\item{[VP]}
Van Proeyen, L., Local zeta functions for ideals and the monodromy
conjecture, thesis 2008 Katholieke Universiteit Leuven.

\item{[VV]}
Van Proeyen, L.\ and Veys, W., The monodromy conjecture for zeta
functions associated to ideals in dimension two,
Annales Inst.\ Fourier, 60 (2010), 1347--1362.

\item{[Ve1]}
Veys, W., On the poles of Igusa's local zeta function for curves,
J.\ London Math.\ Soc.\ (2) 41 (1990), 27--32.

\item{[Ve2]}
Veys, W., Poles of Igusa's local zeta function and monodromy,
Bull.\ Soc.\ Math.\ France 121 (1993), 545--598.

\item{[Ve3]}
Veys, W., Determination of the poles of the topological zeta
function for curves, Manuscripta Math.\ 87 (1995), 435--448.

\item{[Ve4]}
Veys, W., Vanishing of principal value integrals on surfaces,
J.\ Reine Angew.\ Math.\ 598 (2006), 139--158.

\item{[Yu]}
Yuzvinsky, S., Cohomology of the Brieskorn-Orlik-Solomon algebras,
Comm.\ Algebra 23 (1995), 5339--5354.

{\sfont\baselineskip=10pt
\ms
Department of Mathematics, University of Notre Dame, IN 46556,
USA

\sk
RIMS Kyoto University, Kyoto 606-8502 Japan

\sk
Department of Mathematics, University of Oregon, Eugene, Oregon
94703, USA

\ms\vers}}
\bye